\documentclass{amsart}
\usepackage{amsmath,amsthm,amssymb,amsfonts}
\usepackage[colorlinks=true]{hyperref}

\hypersetup{urlcolor=blue, citecolor=blue}

\def\arXiv#1{   {\href{http://arxiv.org/abs/#1}
   {{\mdseries\ttfamily arXiv:#1}}}}

\let\MR\mr

\def\doi#1{   {\href{http://dx.doi.org/#1}
   {{\mdseries\ttfamily DOI}}}}

\newcommand{\al}{\alpha}    \newcommand{\be}{\beta}
\newcommand{\de}{\delta}    \newcommand{\De}{\Delta}
\newcommand{\vep}{\epsilon}  \newcommand{\ep}{\varepsilon}

\newcommand{\R}{\mathbb{R}}

\newcommand{\ti}{\tilde }

\newcommand{\pt}{\partial_t}\newcommand{\pa}{\partial}
\newcommand{\les}{{\lesssim}}
\newcommand{\beeq}{\begin{equation}}\newcommand{\eneq}{\end{equation}}

\def\CO{\mathcal {O}}

\newenvironment{prf}{\noindent {\bf Proof.} }{\endprf\par}
\def \endprf{\hfill  {\vrule height6pt width6pt depth0pt}\medskip}

\numberwithin{equation}{section}

\def\<{\langle}             \def\>{\rangle}
\def\({\left(}                 \def\){\right)}

\newtheorem{theorem}{Theorem}[section]
\newtheorem{corollary}[theorem]{Corollary}

\newtheorem{lemma}[theorem]{Lemma}

\theoremstyle{definition}

\newtheorem{remark}[theorem]{Remark}

\theoremstyle{definition}
\title[Global existence for quasilinear wave equations]
      {Global existence for some 4-D quasilinear wave equations with low regularity}

\author{Mengyun Liu}
\address{School of Mathematical Sciences\\                Zhejiang University\\                Hangzhou 310027, P. R. China}
\email{mengyunliu@zju.edu.cn}

\author{Chengbo Wang}\address{School of Mathematical Sciences\\                Zhejiang University\\                Hangzhou 310027, P. R. China}\email{wangcbo@zju.edu.cn}
\urladdr{http://www.math.zju.edu.cn/wang}

\thanks{The second author was supported in part by National Support Program for Young Top-Notch Talents.}
\date{\today}
\dedicatory{} \commby{}

\begin{document}

\begin{abstract}
In this paper, we prove the global existence for some 4-D quasilinear wave equations with small, radial data in $H^{3}\times H^{2}$. The main idea is to exploit local energy estimates with variable coefficients, together with the trace estimates.
\end{abstract}



\maketitle
%\tableofcontents

%%% Section 1 %%%
\section{Introduction}
In this paper, we are interested in the global solvability of small-amplitude spherically symmetric solutions for the Cauchy problem of the quasilinear wave  equations with low regularity
\begin{align}
\label{ea1}
&\partial_t^2 \phi - \Delta \phi + h(\phi)\Delta \phi = F(\partial \phi)
\quad  (t,x)\in (0,\infty)\times {\mathbb R}^4\\
\label{ea2}
&\phi(0,x)=f(x) \in H^{s}_{{\rm rad}}({\mathbb R}^4), \partial_{t}\phi(0, x)=g(x) \in H^{s-1}_{{\rm rad}}({\mathbb R}^4)
\end{align}
where $h$ is smooth and $h(0) = 0$. We also assume that, for some even number $p\ge 2$, the nonlinearity $F$ takes the form
$$F(\partial \phi) = a |\partial_{t}\phi |^{p} + b|\nabla \phi|^{p}, p\in 2\mathbb{N}, p\ge 2$$
for some $(a, b)\in \R^2$. Here,  $H^{s}_{{\rm rad}}$ stands for the space of spherically symmetric functions lying in the usual Sobolev space $H^s$.

Before presenting our main result, let us first give a brief review of the history, in a broader context.
When $h=0$, the equation (\ref{ea1}) is semiliner:
\begin{equation}
\label{lm1}
\phi_{tt}-\Delta\phi=a|\partial_{t}\phi|^{p}+ b|\nabla \phi|^{p},
\end{equation} which is scale-invariant in the homogeneous Sobolev space $\dot H^{s_c}$ with $$s_{c}=\frac{n}{2}+1-\frac{1}{p-1}\ .$$
With general $p\in \R$ and $p>1$, the local theory has been well understood. It is known that, with general spatial dimension $n$,
the problem is locally well-posed in $H^{s}$ only if $s \ge s_{c}$ (see e.g.  Fang-Wang \cite{Fang2005Local}).
On the other hand, the problem is locally well-posed in $H^{s}$, for $s > \max(s_{c}, \frac{n+5}{4})$, at least when $p\ge 3$ or $p=2$ (see Tataru \cite{Tataru1999On}, Fang-Wang \cite{Fang2005Local} and references therein).
In addition, if $p \geq 3$ with $n\ge 2$ or $p=2$ with $n \geq 3$, the problem is locally well-posed in $H^{s}$ for $s > s_{c}$ when the initial data have radial symmetry or certain amount of angular regularity (see Fang-Wang \cite{Fang2010} ($p\ge 3$), Sterbenz \cite{Sterbenz2007Global} ($n\ge 4$), Hidano-Jiang-Lee-Wang \cite{HJLW} ($n=3$ and $p=2$)).

Concerning the global solvability of \eqref{lm1} with classical small data, it is related to the Glassey conjecture and it is known that
the problem admits global small solutions in general, only if $p > 1+\frac{2}{n-1}$ (see Zhou \cite{Zhou2001BLOW} and references therein). The sufficiency of the condition was known for $n=2, 3$ (Hidano-Tsutaya \cite{Hidano1995Global}, Tzvetkov \cite{Tzvetkov1998Existence}) and $n\ge 4$ in the radial case
(Hidano, Wang and Yokoyama \cite{Hidano2011The}).
The problem of long time existence with low regularity for \eqref{lm1} has been well investigated, see
\cite{Hidano2006Space}, \cite{Sterbenz2007Global},  \cite{Fang2010},  \cite{HJLW} and references therein.

Turning to the general quasilinear wave equation \eqref{ea1}, the local theory was well developed for the case $p=2$ and local well-posedness in $H^s$ is known when $s > \max(\frac{n+1}{2}, \frac{n+5}{4})$ for $n\le 5$, and
$s> \frac{n+1}{2}+\frac{1}{6}$ for $n\ge 6$, see Smith-Tataru \cite{Smith2005Sharp} and references therein.
As the nonlinearity is algebraic, the problem admits global small solutions provided that  $n\ge 4$ or $n=3$ with $F=0$, see e.g. \cite{Ho91}, \cite{L3} and references therein.

There is not much work on the long time existence with low regularity for the quasilinear equation \eqref{ea1}. The only known results are for the case $n=3$ and $p=2$ with radial data. When $F=0$, Zhou and Lei \cite{Zhou2007Global} proved global existence for compactly supported
$H^2_{\rm rad}\times H^1_{\rm rad}$ data.
For general $F$ with $p=2$, the problem was investigated by Hidano, Wang and Yokoyama \cite{Hidano2012On}, where almost global existence up to $\exp(c/\ep)$ was shown for $H^{2}_{{\rm rad}}\times H^{1}_{{\rm rad}}$ data of size $\ep$. When $n$ is odd and $n\ge 5$, the high dimensional global existence with small $H^{(n+1)/2}_{{\rm rad}}\times H^{(n-1)/2}_{{\rm rad}}$ data has also been proven in \cite{HWY-10un}. In this paper, when $n=4$ and the initial data are radial, we would like to investigate the corresponding problem with minimal regularity.

Now let us state our main result.
\begin{theorem}
\label{ta1}
Let $p\in \R$, $p$ is even and $p\ge 2$. Considering \eqref{ea1}-\eqref{ea2} with $s=3$. There exist a constant
 $\ep_0>0$ such that
 for any data $(f,g)$ with
\begin{equation*}
 \|(\nabla f,g )\|_{H^2({\mathbb R}^4)} \le \ep_0\ ,
\end{equation*}
 the problem has a unique solution $\phi\in L^\infty([0,\infty); H^{3}_{{\rm rad}})\cap C^{0,1}([0,\infty); H^{2}_{{\rm rad}})$ with $r^{-1/4}(\pa, \frac{1}{r})\phi \in L^2([0,T]; L^2(\R^4))$ for any $T\in (0,\infty)$. Here $r=|x|$, $\pa=(\pt,\nabla)$ denotes the space-time derivatives, and we denote by $C^{0,1}$ the space of Lipschitz continuous functions.
% \cap C([0,\infty); H^{2}_{{\rm rad}})\cap C^1([0,\infty); H^{1}_{{\rm rad}})$.
\end{theorem}
\begin{remark}
Following the argument in \cite[Section 4]{Hidano2012On}, the solution could be shown to be in the even smaller space $C([0,\infty); H^{3}_{{\rm rad}})\cap C^1([0,\infty); H^{2}_{{\rm rad}})$. We leave the details to the interested reader.
\end{remark}

\begin{remark}
The regularity assumption on the initial data is sharp in general,
as $\sup_{p\ge 2} s_c=3$.
\end{remark}

\begin{remark}
As we can see, our proof is robust enough to apply for situations when the coefficients $a$, $b$ replaced by sufficiently regular (say $C^3$) functions of $\phi$, and $p$ can be any real number $p\ge 2$. Finite combinations of nonlinearities with different $p$'s could also be allowed.
\end{remark}

Let us conclude the introduction by describing the strategy of the proof. We basically follow the approach that appeared in  \cite{Hidano2012On} to give the proof. We mainly exploit a variant of the local energy estimates for wave equations with variable coefficients, due to \cite{Metcalfe2006Long} and \cite{Hidano2012On}. Such kind of estimates have been proven to be extremely useful for the wave equations, see e.g., \cite{MTT, LMSTW, W15} for more history and applications of such estimates.
To effectively exploit the spherically symmetric assumption, trace estimates (Lemma \ref{ld2}, Corollary \ref{th4}) play an important role, to gain spatial decay.

\subsubsection*{Notations} We close this section by listing the notation.
Let $S_T = (0,T) \times {\mathbb R}^4$,
 we set
$X_{m}(T)= \bigcap_{j=0}^{1} C^{j}([0,T);H_{{\rm rad}}^{m-j}({\mathbb{R}}^{4})), m=1,2,3.$
For each $\phi \in X_m(T)$, we let
$ \|\phi\|_{X_m(T)} =
\|\phi\|_{L^{\infty}([0,T); L^2({\mathbb R}^4))} + \|\phi\|_{E_m(T)}$,
where
$$
 \|\phi\|_{E_m(T)} = \|\partial \nabla^{\le m-1} \phi\|_{L^{\infty}((0,T); L^2({\mathbb R}^4))}:=
\sum_{|\al|\le m-1} \|\partial \nabla^{\al} \phi\|_{L^{\infty}((0,T); L^2({\mathbb R}^4))}.
$$
 Secondly, we define
$Y_m(T) = \{\phi\,;\, \phi \in L^1_{\rm loc}(S_T),\enskip
 \|\phi\|_{Y_m(T)} < \infty\}$ for $m\ge 1$, where
$$ \|\phi\|_{Y_m(T)} = (1+T)^{-1/4} \|r^{-1/4}\tilde\pa \nabla^{\le m-1}\phi\|_{L^2(S_T)},
$$ where
$\tilde \pa \phi=(\pa \phi, \frac{1}{r}\phi)$.
Lastly, with $\<r\>=\sqrt{1+r^2}$, we set
%$Z_m(T) = \{\phi\,;\, \phi \in L^1_{\rm loc}(S_T),\enskip \|\phi\|_{Z_m(T)} < \infty\}$ for $m\ge 1$, where for some $\de\in (0, 1/4)$ to be determined,
$$
 \|\phi\|_{Z_m(T)}=
 \|r^{-1/4} \langle r \rangle^{-1/4-\delta}\tilde\pa
\nabla^{\le m-1}
  \phi\|_{L^2(S_T)}\ ,  \de\in (0, 1/4)\ .$$
When $T=\infty$, it is understood that $E_m=E_m(\infty)$, $Z_m=Z_m(\infty)$. We will also use $A\les B$ to stand for $A\le C B$ where the constant $C$ may change from line to line.

% for some $\de\in (0, 1/4)$ to be determined,
%%% Section 2 %%%

\section{Preliminaries}

In this section, we collect various basic estimates to be used.
We will use the following version of the weighted Sobolev estimates, which essentially are consequences of the well-known trace estimates, see, e.g., \cite[(1.3), (1.7)]{FaWa} and references therein.
\begin{lemma}[Trace estimates]
 \label{ld2}
Let $n\geq 2 $ and $1/2 < s < n/2 $, then
\begin{equation}
\| r^{n/2-s}u\|_{L_{r}^{\infty}L_{\omega}^{2}} \les\|u\|_{\dot{H}^{s}},\|r^{(n-1)/2}u\|_{L_{r}^{\infty}L_{\omega}^{2}}\les \|u\|_{H^{s}}\ .
\label{ed23}
\end{equation}
In particular, we have
\begin{equation}
\|r^{n/2 -s}\langle r\rangle^{s-1/2}u\|_{L_{r}^{\infty} L^2_\omega} \les\|u\|_{H^{s}}.
\label{ed22}
\end{equation}
\end{lemma}
\begin{corollary}
\label{th4}
Let $n\ge 4$, $\alpha$ be multi-indices with $|\alpha|=1, 2$, and $1/2<s< n/2+1-|\alpha|$,  we have for any radial functions $u$,
\begin{equation}
|\nabla^{\alpha}u|\leq C r^{s-\frac{n}{2}}\langle r \rangle^{\frac{1}{2}-s}
\sum_{|\be|=|\al|}
\|\nabla^{\be}u\|_{H^{s}}.
\label{edx915}
\end{equation}
\end{corollary}
\begin{prf}
Since u is radial, we have $\nabla u=\frac{x}{r}\partial_{r}u$, then
\begin{equation*}
|\pa_j u|=|\frac{x_j}{r}\partial_{r}u|\le |\partial_{r}u|=A_{n-1}^{-\frac{1}{2}}\|\partial_{r}u\|_{L_{\omega}^{2}}=
A_{n-1}^{-\frac{1}{2}}\|\nabla u\|_{L_{\omega}^{2}},
\end{equation*}
with $A_{n-1}=|{\mathbb{S}}^{n-1}|$, by (\ref{ed22}),
\begin{equation}
\|r^{n/2-s}\<r\>^{s-1/2}\pa_j u\|_{L_{x}^{\infty}} \les \|r^{n/2-s}\<r\>^{s-1/2}\nabla u\|_{L_{r}^{\infty}L_{\omega}^{2}}
\les \|\nabla u\|_{H^{s}}, 1/2<s<n/2.
\label{edx916}
\end{equation}
Which is \eqref{edx915} for $|\alpha|=1$.
For $|\alpha|=2$, with $\omega_{i}=\frac{x_{i}}{r}$, we notice that $u_{x_{i}x_{j}}=\omega_{i}\omega_{j}u_{rr}+\frac{\delta_{ij}-\omega_{i}\omega_{j}}{r}u_{r}$. Based on this fact, we can get
$$
|u_{x_{i}x_{j}}|\les |u_{rr}|+\frac{1}{r}|u_{r}|,
|u_{rr}|\les |\nabla ^2 u|+\frac{1}{r}|\nabla u|\ .
$$
Thus if we have $s\in (1/2, (n-2)/2)$ such that $s, s+1\in (1/2, n/2)$, we have
\begin{eqnarray*}
\|r^{n/2-s}\<r\>^{s-1/2} u_{x_{i}x_{j}}\|_{L_{r}^{\infty}L_{\omega}^{\infty}}&\les& \|r^{n/2-s}\<r\>^{s-1/2} (u_{rr}, u_r/r)\|_{L_{r}^{\infty}L_{\omega}^{\infty}}\\
&\les& \|r^{n/2-s}\<r\>^{s-1/2} (u_{rr}, u_r/r)\|_{L_{r}^{\infty}L_{\omega}^{2}}\\
&\les& \|r^{n/2-s}\<r\>^{s-1/2} (\nabla^2 u, \nabla u/r)\|_{L_{r}^{\infty}L_{\omega}^{2}}\\
&\les&  \|\nabla^2 u\|_{{H}^{s}}+\|\nabla u\|_{\dot{H}^{s+1}}
+\|\nabla u\|_{\dot{H}^{3/2}}\\
&\les& \|\nabla^2 u\|_{H^s},
\end{eqnarray*}
which completes the proof.
\end{prf}

We will also need the following weighted Sobolev estimates.
\begin{lemma}
 \label{thm-w-decay}
Let $n=4$ and $m\in \R$, then for radial functions $u$, we have
\begin{equation}
\| \<r\>^{3/4+m}u\|_{L_{x}^{4}} \les\|\<r\>^m \nabla^{\le 1} u\|_{L^{2}}, \| \<r\>^{3/4+m}\nabla^k u\|_{L_{x}^{4}} \les\|\<r\>^m \nabla^{\le 1}\nabla^k u\|_{L^{2}}, k\ge 1\ .
\label{eq-L4decay}
\end{equation}
\end{lemma}
\begin{prf}
The estimates are direct consequences of the particular case \cite[Lemma 2.2 (2.3)]{W15}, which states that for general functions $u$ in $\R^4$,
$$\| \<r\>^{3/4+m}u\|_{L_{x}^{4}} \les \|\<r\>^m \nabla^{\le 1} u\|_{L^{2}}+
\|\<r\>^m \Omega u\|_{L^{2}}\ ,
$$ where $\Omega=\{\Omega_{ij}=x^i\pa_j-x^j\pa_i\}$.
As $u$ is radial, $\Omega u=0$ and so is the first estimate.
For the second estimate with $k=1$, we notice  $\Omega_{ij} \pa_k u=\pa_k\Omega_{ij}  u-\de_{ik}\pa_j u-\de_{jk}\pa_i u
=-\de_{ik}\pa_j u-\de_{jk}\pa_i u
$ and then get
$$\| \<r\>^{3/4+m}\pa_k u\|_{L_{x}^{4}} \les \|\<r\>^m \nabla^{\le 1} \pa_k u\|_{L^{2}}+
\|\<r\>^m \Omega \pa_k u\|_{L^{2}}
\les
 \|\<r\>^m \nabla^{\le 1} \nabla u\|_{L^{2}}
\ .
$$  The same argument gives us the second estimate for general $k$.
\end{prf}

We will also need a variant of the local energy estimates (see \cite[Section 2]{Hidano2012On}).
Let $h^{\alpha \beta} \in C^1([0,T]\times {\mathbb R}^n)$,
$h^{\alpha \beta}=h^{\beta \alpha}$, $0\le \al,\be\le n$ satisfying
 \begin{equation}
 \sum_{\alpha,\beta = 0}^n |h^{\alpha \beta}(t,x)| \leq \frac{1}{2}, \forall (t,x) \in [0,T]\times {\mathbb R}^n.
\label{eb3}
\end{equation}
 Consider the  linear wave equation with variable coefficients
\begin{equation}
 \partial_t^2 \phi - \Delta \phi + \sum_{\al,\be=0}^n h^{\alpha \beta}(t,x)\partial_\alpha
  \partial_\beta \phi = F(t,x)
 \quad
\mbox{in}\, \,(0,T)\times{\mathbb R}^n,
\label{eb1}
\end{equation}
with the initial data
\begin{equation}
 \phi(0,\cdot) = f \in H^1({\mathbb R}^n), \quad
\partial_t \phi (0, \cdot) = g \in L^2({\mathbb R}^n)
.
\label{eb2}
\end{equation}

\begin{lemma}[Local energy estimates]
\label{tb1}
Let $n \geq 3$, $\mu \in (0,1/2)$, $T\ge 1$ and $\phi \in C^2([0,T]\times {\mathbb R}^n)$ be a solution of the initial value problem \eqref{eb1}-\eqref{eb2} satisfying the condition \eqref{eb3}. Suppose also that $\phi \in C^0([0,T];H^1({\mathbb R}^n)) \cap
 C^1([0,T];L^2({\mathbb R}^n))$ and
$ |\tilde\partial \phi||F| \in L^1((0,T)\times {\mathbb R}^n)$.
 Then, for any $\delta > 0$, there exists $C_0>0$, which is independent of $T\ge 1$, such that we have
\begin{eqnarray}
&&
\|\pa \phi\|_{L^\infty([0,T]; L^2)}^2+
\|r^{-1/2 + \mu}\langle r\rangle^{-\mu-\delta} \tilde\partial \phi\|^2_{L^2(S_T)}
+T^{-2\mu}\|r^{-1/2 + \mu} \tilde\partial \phi\|^2_{L^2(S_T)}
\label{eb5} \\
& \leq & C_0 \|(\nabla f, g)\|^2_{L^2({\mathbb R}^n)}
+ C_0 \int_0^T\int_{{\mathbb R}^n}
|\tilde\partial \phi|\left(|F|
+|\pa \phi| \left(|\partial h|+ \frac{| h|}{r^{1-2\mu}\langle r\rangle^{2\mu}}
\right)\right) dxdt.\nonumber
 \end{eqnarray}
%Here $C_1$ is a positive constant depending only on $n$, $\mu$, and $\delta$.
\end{lemma}

\section{Global Existence}
In this section, we give the proof of the global existence for the quasilinear wave equation \eqref{ea1}-\eqref{ea2}, by exploiting the local energy estimates, as well as the spatial decay provided by the trace estimates.
For this purpose, our first step is to define a sequence of approximate solutions as follows.

\subsection{Approximate solutions}
Firstly, we choose a radially symmetric function
$\rho \in C_0^{\infty}({\mathbb R}^4)$ with $\int_{{\mathbb R}^4} \rho(x)dx = 1$, and set
$\rho_{k} (x) = 2^{4k} \rho (2^{k}x)$.
Based on $\rho$, we define standard sequence of $C^\infty$,  radially symmetric, approximate functions to $(f,g)$,
\begin{equation}
 f_k(x) = \rho_{k} * f(x),
g_k(x) = \rho_{k} * g(x), k\ge 0\ .
\label{ed11}
\end{equation}
As is well known,
$\|f_k - f\|_{H^3} \rightarrow 0$,
$\|g_k - g\|_{H^2} \rightarrow 0$,
as $k \to \infty$.
We also have
\begin{equation}
 \sum_{k = 1}^{\infty} \bigl( \|\nabla f_k
 - \nabla f_{k-1}\|_{L^2({\mathbb R}^4)}
+ \|g_k
 - g_{k-1}\|_{L^2({\mathbb R}^4)}\bigr)
< \infty.
\label{ed12}
\end{equation}
Indeed, we can easily check this property by using
$\|\rho_{k} * \varphi - \varphi\|_{L^2} \leq C2^{-k}\|\varphi\|_{H^1}$.

With $(f_k,g_k)$ as data, we use a standard iteration to define the sequence of approximative solutions. Let $\phi_{-1} \equiv 0$ and define $\phi_k$ ($k\ge 0$) recursively by solving
\beeq
\left\{\begin{array}{l}\partial_t^2 \phi_k - \Delta \phi_k + h(\phi_{k-1})
  \Delta \phi_k =  F(\partial \phi_{k-1}), (t,x)\in S_T,
 \\
\phi_k(0,\cdot) = f_k,\quad
 \partial_t \phi_k(0,\cdot) = g_k.\end{array}
 \right.
 \label{ed14}% \label{ed15}
\eneq
In order to ensure that the sequence $\{\phi_k\}$ is well-defined,
we will have to assume that $\varepsilon = \|\nabla f\|_{H^2} + \|g\|_{H^2}$ is small enough.
Applying a standard existence, uniqueness and regularity theorem,
 we will see that, for all $k = 0,1,2,\ldots$,
$\phi_k$ is well defined,
radially symmetric for all time, and satisfies
\begin{equation}
\phi_k \in C^{\infty}(\overline{S_T}) \cap X_3(T)\cap
 Y_3(T) \cap
Z_3(T) .
  \label{ed16}
\end{equation}

\subsection{A-priori estimates} To obtain the desired properties of the approximate solutions, we give some a-priori estimates for the equations.

Let $n = 4$ and $\mu = 1/4$ in Lemma \ref{tb1}, then
\eqref{eb5} reads as follows,
\beeq
\label{lm2}
\|\phi\|^2_{{E_1(T)\cap Y_1(T) \cap Z_1(T)}}
\les \|(\nabla f,g)\|^2_{L^2({\mathbb R}^4)} + \int_0^T \int_{{\mathbb R}^4}
|\tilde \partial \phi|
\left(|F| + |\pa \phi| |\bar \partial h|\right) dxdt.
\eneq
where we have introduced $\bar\pa h=(\pa h,  \frac{ h}{r^{1/2}\langle r\rangle^{1/2}})$
and $\|\phi\|_{E_1(T)\cap Y_1(T) \cap Z_1(T)}=
\|\phi\|_{E_1(T)} +\|\phi\|_{Y_1(T)} + \|\phi\|_{Z_1(T)}$
 for brevity.

By \eqref{lm2} and the definitions of $Y_1$, $Z_1$ norms, we immediately get the following lemma.
%%% Lemma ld7 %%%
\begin{lemma}
 \label{ld7}
Let $T\ge 1$, $\phi \in C^2(\overline{S_T}) \cap X_1(T)$ be a solution of the wave equation
\begin{equation}
 \label{edx9}
\partial_t^2 \phi - \Delta \phi + h(t,x)\Delta \phi = F
\quad \mbox{in}\,\, S_T,
\end{equation}
with initial data $(f,g) \in H^1({\mathbb R}^4)\times L^2({\mathbb R}^4)$.
Assume $h \in C^1(\overline{S_T})$ satisfying
\begin{equation}\label{eq-cond-h}
 \|h\|_{L^{\infty}(S_T)} \leq \frac{1}{8},\ \| r^{1/2}\langle r \rangle^{1/2+2\de}\bar\partial h\|_{L^{\infty}(S_T)} < \infty .
\end{equation}
Then we have
\begin{eqnarray}
\label{edx12a}
&&\|\phi\|_{E_1(T)\cap Y_1(T) \cap Z_1(T)}\\
 &\les& \|(\nabla f,g)\|_{L^2(\mathbb{R}^{4})}+T^{1/4}\|r^{1/4}F\|_{L^2(S_T)} + T^{1/2}\|\phi\|_{Y_1(T)}
\| r^{1/2}\bar\partial  h\|_{L^{\infty}(S_T)}\ ,
\nonumber
\end{eqnarray}
\begin{eqnarray}
\label{edx12b}
&&\|\phi\|_{E_1(T)\cap Y_1(T) \cap Z_1(T)}\\
 &\les& \|(\nabla f, g)\|_{L^2(\mathbb{R}^{4})} +
  \|r^{1/4}\langle r \rangle^{1/4+\delta}F\|_{L^2(S_T)}  +\|\phi\|_{Z_1(T)}\| r^{1/2}\langle r \rangle^{1/2+2\de}\bar\partial h\|_{L^{\infty}(S_T)}\ .
 \nonumber
\end{eqnarray}
\end{lemma}
\begin{prf} These estimates are direct consequences of \eqref{lm2} by the Cauchy-Schwartz inequality. For illustration, let us present the proof for
\eqref{edx12a}. At first, we observe that the second term on the right hand side of \eqref{lm2} could be controlled by
\begin{eqnarray*}
%\int_0^T \int_{{\mathbb R}^4}|\tilde \partial \phi|\left(|F| + |\pa \phi| |\bar \partial h|\right) dxdt
&& \int_0^T \int_{{\mathbb R}^4}
r^{-1/4}|\tilde \partial \phi|
r^{1/4}\left(|F| +|\pa \phi| |\bar \partial h|\right) dxdt\\
&\les&
\|r^{-1/4}\tilde \partial \phi\|_{L^2(S_T)}
(\|r^{1/4}F\|_{L^2(S_T)}+\|r^{1/4} \pa \phi   \bar \partial h\|_{L^2(S_T)})\\
&\les&
T^{1/4}\|\phi\|_{Y_1(T)}
\left(\|r^{1/4}F\|_{L^2(S_T)}+T^{1/4}\|\phi\|_{Y_1(T)}\|r^{1/2}  \bar \partial h\|_{L^\infty(S_T)}\right)\ .
\end{eqnarray*}
Plugging this estimate to \eqref{lm2}, we arrive at
\eqref{edx12a} by absorbing the term $\|\phi\|_{Y_1(T)}$.
This completes the proof.
\end{prf}

With help of Lemma \ref{ld7}, we could prove the following higher order estimate.
\begin{lemma}[Higher order estimate]
 \label{ld6}
Let $\phi, {\tilde \phi} \in
C^{\infty}(\overline{S_T}) \cap X_3(T)\cap Z_3(T)$.
Assume that they satisfy
\begin{equation}
 \label{edx1}
 \partial_t^2 \phi - \Delta \phi + h({\tilde \phi})\Delta \phi
= F(\partial {\tilde \phi})
\ , \, (t,x)\in S_T\ .
\end{equation}
Then there exists $C_1>0$, independent of $T$, such that we have
\beeq
 \label{edx2}
\|\phi\|_{E_3(T) \cap Y_3(T)\cap Z_3(T)}\les \|\partial \phi(0,\cdot)\|_{H^2} +
 %\|\partial {\tilde \phi}(0,\cdot)\|^2_{H^2 }+
   \|{\tilde \phi}\|_{E_3(T)}(\|{\tilde \phi}\|_{Z_3(T)}+\|{\phi}\|_{Z_3(T)})
\eneq
provided that
$\|{\tilde \phi}\|_{E_3(T)} \leq C_1$.
\end{lemma}
\begin{prf}
Since $W^{2,4}({\mathbb{R}}^{4})\hookrightarrow L^{\infty}({\mathbb{R}}^{4})$,  ${\dot{H}}^{1}({\mathbb{R}}^{4})\hookrightarrow L^{4}({\mathbb{R}}^{4})$, we have
$$\|u\|_{L^\infty}\les\|\nabla^{\le 2}u\|_{L^{4}} \les \|\nabla^{\le 2}u\|_{\dot H^{1}}\ .$$
As $h\in C^{3}$ and $h(0)=0$, there exists a positive constant $C_{1}$ such that $|h(\tilde \phi)|\leq \frac{1}{8}$ in $S_T$,
provided that $\|\tilde \phi\|_{E_3(T)}\leq C_{1}$. Also, we have $h^{j}(\ti \phi)\les 1$ for $j\le 3$.
By the inequality \eqref{ed23}, as $\de\in (0, 1/4)$, we see that
$$\|r^{2\delta}h(\tilde \phi)\|_{L^{\infty}(S_{T})}\les \|h(\tilde{\phi})\|_{L^{\infty} \dot{H}^{2-2\de}}\les \|h(\tilde{\phi})\|_{L^{\infty} (\dot{H}^{1}\bigcap \dot{H}^{3})}\les \|\tilde \phi\|_{E_{3}(T)}<\infty,$$
$$
\|r^{\frac{1}{2}}\langle r\rangle^{\frac{1}{2}+2\delta}\partial h(\tilde \phi)\|_{L^{\infty}(S_{T})}\les\|\partial h(\tilde{\phi})\|_{L^{\infty} H^{\frac{3}{2}}}\les
\|\partial h(\tilde{\phi})\|_{L^{\infty} H^{{2}}}
\les\|\tilde \phi\|_{E_{3}(T)}<\infty,
$$
which, combined with
$|h(\tilde \phi)|\leq \frac{1}{8}$, gives us
\eqref{eq-cond-h}.
By \eqref{edx915}, we get
\beeq\label{eq-dphi} |\partial\tilde{\phi}(t,\cdot)| \les   r^{-\vep}\langle r\rangle^{\vep-\frac{3}{2}}\|\partial\tilde{\phi}(t,\cdot)\|_{H^{2-\ep}}\les  r^{-\vep}\langle r\rangle^{\vep-\frac{3}{2}}\|\partial\tilde{\phi}(t,\cdot)\|_{H^{2}}, \vep\in (0, 3/2)\ .\eneq
Let $\vep=\frac{1}{2(p-1)}$,
 as $\de\in (0,1/4)$ and $p\ge 2$, we have
$-\vep(p-1)+\frac{1}{4}\geq -\frac{1}{4}$,  $-\frac{3}{2}(p-1)+\frac{1}{2}+\delta \leq -\frac{1}{2}-\delta$.
Since
$|F(\tilde \phi)| \les |\partial \tilde{\phi}|^{p}$,
then
\begin{align*}
\|r^{\frac{1}{4}}\langle r\rangle^{\frac{1}{4}+\delta}F(\tilde \phi)\|_{L^{2}(S_{T})}
&\les  \|r^{-\vep(p-1)+\frac{1}4}\langle r\rangle^{(\vep-\frac{3}{2})(p-1)+\frac{1}{4}+\delta}\partial \tilde{\phi}\|_{L^{2}(S_{T})}\|\tilde \phi\|_{E_{3}(T)}^{p-1}\\
&\les \|\tilde \phi\|_{E_{3}(T)}^{p-1}\|\tilde \phi\|_{Z_{1}(T)}\\
&\les\|\tilde \phi\|_{E_{3}(T)}\|\tilde \phi\|_{Z_{1}(T)}.
\end{align*}
Applying (\ref {edx12b}) to (\ref {edx1}), we obtain
\beeq
\label{ed222}
\|\phi\|_{E_1(T)\cap Y_1(T) \cap Z_1(T)}
%\|\phi\|_{E_1(T)} + \|\phi\|_{Y_1(T)} + \|\phi\|_{Z_1(T)}
\les  \|\partial\phi(0,\cdot)\|_{L^{2}}+\|\tilde \phi\|_{E_{3}(T)}
(\|\tilde \phi\|_{Z_{1}(T)}+\|\phi\|_{Z_{1}(T)}). \eneq

Turning to the
estimate of higher order, we apply $\nabla^{\le 2}$ to the equation to get
\begin{align}
\label{edx11}
&\partial_{t}^{2}(\nabla^{\leq 2}\phi)-\Delta(\nabla^{\leq 2}\phi)+h(\tilde \phi)\Delta(\nabla^{\leq 2}\phi)\\
=&\nabla^{\leq  2}F(\partial \tilde \phi)+\CO(
|h'(\tilde{\phi})\nabla \ti \phi \De \nabla^{\le 1} \phi|+
|h'(\tilde{\phi})(\nabla^2 \ti \phi) \De \phi|+
|h''(\tilde{\phi})(\nabla \ti \phi)^2 \De \phi|)\nonumber\\
=&\nabla^{\leq  2}F(\partial \tilde \phi)+\CO(
|\nabla \ti \phi \De \nabla^{\le 1} \phi|+
|(\nabla^2 \ti \phi) \De \phi|+
|(\nabla \ti \phi)^2 \De \phi|)=R  \nonumber
\end{align}
Then, by Lemma \ref{ld7}, the proof of \eqref{edx2} is reduced to the following estimate
\begin{equation}
\label{ed333}
\|r^{\frac{1}{4}}\langle r\rangle^{\frac{1}{4}+\delta}R\|_{L^{2}(S_{T})}\les \|\tilde\phi\|_{E_{3}(T)}
(\|\tilde\phi\|_{Z_{3}(T)}+\|\phi\|_{Z_{3}(T)}).
\end{equation}

We deal with each item of $R$ separately. At first, for
$\nabla^{\leq  2}F(\partial \tilde \phi)$,
we have
$$|\nabla^{\leq 2} F(\partial \tilde \phi)|\les
|\partial\tilde \phi|^{p-1}|\partial\nabla^{\leq 2} \tilde{\phi}|
+|\pa \ti {\phi}|^{p-2} | \pa \nabla \tilde{\phi}|^2\ .
$$
As $\de \in (0, 1/4)$ and $p\ge 2$,
by choosing $\vep\in (0, \frac{1}{2(p-1)}]$ such that $\vep(p-2)\le 1/4$, we have
$-\vep(p-1)+\frac{1}{4}\geq -\frac{1}{4}$, $-\frac{3}{2}(p-1)+\frac{1}{2}+\delta \leq -\frac{1}{2}-\delta$ and
$-\frac{3}{2}(p-2)-1+\de\le -1/2-\de$.
%by Lemma \ref{ld2} (\ref{ed22}) and Corollary \ref{th4}\begin{equation}\label{lm3}|\partial \nabla\tilde \phi (t,\cdot)|\les r^{-1}\langle r\rangle^{-1/2}\|\partial\nabla \tilde \phi(t,\cdot)\|_{H_{{\rm rad}}^{1}({\mathbb{R}}^{4})},\end{equation}
Then \eqref{eq-dphi} and \eqref{eq-L4decay} gives us
\begin{eqnarray*}
&&\|r^{1/4}\<r\>^{1/4+\de} \nabla^{\leq 2} F(\partial \tilde \phi)\|_{L^2(S_T)}
\\& \les &
\|r^{1/4}\<r\>^{1/4+\de}
|\partial\tilde \phi|^{p-1}|\partial\nabla^{\leq 2} \tilde{\phi}|
\|_{L^2(S_T)}
  +\|r^{1/4}\<r\>^{1/4+\de}
|\pa \ti {\phi}|^{p-2} | \pa \nabla \tilde{\phi}|^2
\|_{L^2(S_T)}\\
& \les &
 \|r^{-\vep(p-1)+\frac{1}4}\langle r\rangle^{(\vep-\frac{3}{2})(p-1)+\frac{1}{4}+\delta}\partial\nabla^{\leq 2} \tilde{\phi}\|_{L^{2}(S_{T})}\|\tilde \phi\|_{E_{3}(T)}^{p-1}
 \\
 &  &
 + \|r^{-\vep(p-2)+\frac{1}4}\langle r\rangle^{(\vep-\frac{3}{2})(p-2)-\frac{1}{2}+\delta}\partial\nabla \tilde{\phi}\|_{L^2_t L_x^{4}(S_{T})}
\|\<r\>^{3/4} \partial\nabla \tilde{\phi}\|_{L^\infty_t L^4_x(S_T)}
\|\tilde \phi\|_{E_{3}(T)}^{p-2}\\
 & \les &
 \|r^{-\frac{1}4}\langle r\rangle^{-\frac{1}{4}-\delta}\partial\nabla^{\leq 2} \tilde{\phi}\|_{L^{2}(S_{T})}\|\tilde \phi\|_{E_{3}(T)}^{p-1}
 \\
 &  &
 + \|\langle r\rangle^{-\frac{3}{2}(p-2)-\frac{1}{4}+\delta}\partial\nabla \tilde{\phi}\|_{L^2_t L_x^{4}(S_{T})}
\| \partial\nabla^{\le 2} \tilde{\phi}\|_{L^\infty_t L^2_x(S_T)}
\|\tilde \phi\|_{E_{3}(T)}^{p-2}
\\
 & \les &
 \|\tilde{\phi}\|_{Z_3(T)}\|\tilde \phi\|_{E_{3}(T)}^{p-1}
 + \|\langle r\rangle^{-\frac{3}{2}(p-2)-1+\delta}\partial\nabla^{\le 2} \tilde{\phi}\|_{L^2_t L_x^{2}(S_{T})}
\|\tilde \phi\|_{E_{3}(T)}^{p-1}\\
 & \les &
 \|\tilde{\phi}\|_{Z_3(T)}\|\tilde \phi\|_{E_{3}(T)}^{p-1}\ .
\end{eqnarray*}

For the second term,
$\nabla\tilde \phi\nabla^{\leq1}\Delta \phi$,
%$\|r^{\frac{1}{4}}\langle r\rangle^{\frac{1}{4}+\delta}h'(\tilde \phi)\partial\tilde \phi\nabla^{\leq1}\Delta \phi\|_{L^{2}(S_{T})}$:
by \eqref{eq-dphi} with $\vep=1/2$, we have
%$$|\partial\tilde \phi|\les ~r^{-\frac{1}{2}}\langle r\rangle^{-1}\|\partial\tilde \phi\|_{H^{\frac{3}{2}}_{{\rm rad}}({\mathbb{R}}^{4})},$$
\begin{align*}
\|r^{\frac{1}{4}}\langle r\rangle^{\frac{1}{4}+\delta}\nabla\tilde \phi \nabla^{\leq1}\Delta \phi\|_{L^{2}(S_{T})}&\les ~\|r^{-\frac{1}{4}}\langle r\rangle^{-\frac{3}{4}+\delta}\nabla^{\leq1}\Delta\phi\|_{L^{2}(S_{T})}\|\partial\tilde \phi\|_{L_{t}^{\infty}H^{\frac{3}{2}}_{{\rm rad}}({\mathbb{R}}^{4})}\\
&\les ~\| \phi\|_{Z_{3}(T)}\|\tilde \phi\|_{E_{3}(T)}.
\end{align*}

Turning to the third term, $|(\nabla^2 \ti \phi) \De \phi|$,
we follow similar argument for $| \pa \nabla \tilde{\phi}|^2$, i.e., using
\eqref{eq-L4decay}, to conclude
 \begin{eqnarray*}
\|r^{\frac{1}{4}}\langle r\rangle^{\frac{1}{4}+\delta}
(\nabla^2 \ti \phi) \De \phi\|_{L^{2}(S_{T})}&\les &
\|\langle r\rangle^{-1/4+\delta}
 \De \phi\|_{L^2_t L^{4}_x(S_{T})}
\|\<r\>^{3/4} \nabla^2 \ti \phi\|_{L^\infty_t L^4_x(S_T)}\\
&\les&
\|\langle r\rangle^{-1+\delta}
\nabla^{\le 1} \De \phi\|_{L^2_t L^{2}_x(S_{T})}
\|\nabla^2\nabla^{\le 1} \ti \phi\|_{L^\infty_t L^2_x(S_T)}\\
&\les & \| \phi\|_{Z_{3}(T)}\|\tilde \phi\|_{E_{3}(T)}.
\end{eqnarray*}

It remains to control the last term, $(\nabla \ti \phi)^2 \De \phi$, for which
we use \eqref{eq-dphi} with $\vep=1/4$ to get
\begin{eqnarray*}
\|r^{\frac{1}{4}}\langle r\rangle^{\frac{1}{4}+\delta} (\nabla \ti \phi)^2 \De \phi \|_{L^{2}(S_{T})}
&\les& \|r^{-1/4}\langle r\rangle^{-9/4+\delta}
 \De \phi\|_{L^{2}(S_{T})}
\|r^{1/4}\<r\>^{5/4}\nabla \ti \phi\|_{L^\infty(S_T)}^2\\
&\les & \| \phi\|_{Z_{3}(T)}\|\tilde \phi\|^2_{E_{3}(T)}\les
 \| \phi\|_{Z_{3}(T)}\|\tilde \phi\|_{E_{3}(T)},
\end{eqnarray*}
where we have used the boundedness of $\|\tilde \phi\|_{E_{3}(T)}$. This completes the proof.
\end{prf}

\subsection{Boundedness of the approximate solutions in $X_3(T) \cap Y_3(T) \cap Z_{3}(T)$}
In this subsection, we show that the approximate solutions $\phi_k$ are globally well-defined and bounded in $X_3(T) \cap Y_3(T) \cap Z_{3}(T)$ for any $T<\infty$.
%%% Section 3 %%%

\begin{lemma}
\label{exd11111}
Let $(f,g) \in H^3_{{\rm rad}}({\mathbb R}^4) \times H^2_{{\rm rad}}({\mathbb R}^4)$ and $\varepsilon = \|\nabla f\|_{H^2} + \|g\|_{H^2}$.
Then there exist positive constants $\varepsilon_1$, $C_2$ such that
if $\varepsilon \leq \varepsilon_1$,
 the functions $\phi_k\in C^\infty\cap C H^s\cap C^1 H^{s-1}$ for any $s\ge 4$ and enjoy the  uniform bound
 %, defined in \eqref{ed14}, are smooth, have the uniform bound
\beeq\label{eq-bd3.3}
 \|\phi_k\|_{E_3}
 + \|\phi_k\|_{Z_3}
  + \sup_{T\ge 1}\|\phi_k\|_{Y_3(T)}
\leq C_2 \varepsilon\ ,
\forall  k \ge 0\ .
\eneq
%for any $T<\infty$,
\end{lemma}
\begin{prf}
The proof proceeds by induction. At first, letting $\phi=\phi_{0}$, $\tilde{\phi}=0$, we have
$$\partial_{t}^{2}\phi_{0}-\Delta\phi_{0}=0, \phi_0(0)=f_0, \pt\phi_0(0)=g_0\ ,
$$
which is global solvable, well-defined and smooth in $S_T$ for any $T<\infty$. Moreover,
applying Lemma \ref{ld6},
we get
\begin{equation*}
\|\phi_0\|_{E_3(T)} + \|\phi_0\|_{Y_3(T)} + \|\phi_0\|_{Z_3(T)} \leq
C (\|\nabla f_0\|_{H^2} + \|g_0\|_{H^2}),
\end{equation*}
for some $C>0$, independent of $T\ge 1$.
As we know from the definition \eqref{ed11} of $f_k$, $g_k$ that
$(f_k, g_k)\in H^s\times H^{s-1}$ for any $s\ge 1$ and
\begin{equation*}
\|\nabla f_k\|_{H^2} + \|g_k\|_{H^2}
  \leq \|\nabla f\|_{H^2} + \|g\|_{H^2}
 = \varepsilon, \forall k\ge 0\ ,
\end{equation*}
 we have
\begin{equation*}
\|\phi_0\|_{E_3(T)} + \|\phi_0\|_{Y_3(T)} + \|\phi_0\|_{Z_3(T)} \leq
C \varepsilon.
\end{equation*}
In order to define $\phi_{1}$ by \eqref{ed14},
we need to verify that $\phi_0$ satisfies $\|\phi_0\|_{E_3(T)} \leq C_{1}$ so that $|h(\phi_0)|\le 1/8$.
Let $C_2=2C$ and $C_2\varepsilon_1\le C_1$, we make the inductive assumption that for some $m\ge 1$, we have
for any $k\le m-1$,
$\phi_k\in C^\infty\cap C H^s\cap C^1 H^{s-1}$ for any $s\ge 4$
with the bound \eqref{eq-bd3.3} satisfied. With $\phi=\phi_{m}$,
$\tilde{\phi}=\phi_{m-1}$ in Lemma \ref{ld6},
we see that $F(\pa\ti \phi)\in L^1([0,T); H^{s-1})$,
$h(\ti \phi)\in C^\infty$, $|h(\ti \phi)|\le 1/8$
 and $\phi=\phi_{m}$ is
  global solvable, well-defined, smooth in $S_T$ and
  $\phi_k\in C H^s\cap C^1 H^{s-1}$ for any $s\ge 4$.
Moreover, by \eqref{edx2},
 we have
 \begin{align*}
&\|\phi_m\|_{E_3(T)} + \|\phi_m\|_{Y_3(T)} + \|\phi_m\|_{Z_3(T)}\\
\leq& C\|\partial \phi_{m}(0,\cdot)\|_{H^2({\mathbb R}^4)}+C \| \phi_{m-1}\|_{E_3(T)}(\| \phi_{m}\|_{Z_3(T)}+\|{\phi_{m-1}}\|_{Z_3(T)})\\
\leq & C\varepsilon +C (C_{2}\varepsilon)^{2}+C C_{2}\varepsilon\|\phi_m\|_{Z_3(T)}.
\end{align*}
Now if we choose $\varepsilon_{1}>0$ such that
\begin{equation*}
4C^2 \varepsilon_1 \leq \frac{1}{3} \ ,\
2 C\ep_1\le C_1\ ,
\end{equation*}
then for any $\varepsilon\leq \varepsilon_{1}$,
\begin{equation*}
\|\phi_m\|_{E_3(T)} + \|\phi_m\|_{Y_3(T)} + \|\phi_m\|_{Z_3(T)}\leq C_{2}\varepsilon.
\end{equation*}
%As all of the estimates are independent of $T\ge 1$,
This completes the proof by induction.
\end{prf}

\subsection{Convergence in  $X_1(T) \cap Y_1(T) $}
In this subsection, we show that the approximate solutions are convergent to the desired solution of the quasilinear problem. At first, let us give the following lemma.
\begin{lemma}
 \label{ld9}
Let $\phi^{(i)}, {\tilde \phi}^{(i)} \in C^{\infty}({S_T}) \cap X_3(T)$ $(i = 1,2)$.
Assume that they satisty
\begin{equation}
 \partial_t^2 \phi^{(i)} - \Delta \phi^{(i)}
+ h({\tilde \phi}^{(i)})\Delta \phi^{(i)} =
 F(\partial {\tilde \phi^{(i)}})
\quad \mbox{in}\,\,S_T
\nonumber
\end{equation}
and $\|{\tilde \phi}^{(i)}\|_{E_3(T)} \leq C_{1}$ for $i = 1,2$. Then there exists $C_3>0$ such that for any $T\ge 1$, we have
\begin{eqnarray}
&&\|\phi^{(1)} - \phi^{(2)}\|_{E_1(T)\cap Y_1(T)\cap Z_1(T)}  \label{ed45}\\
&\leq& C_3\|\partial \phi^{(1)}(0,\cdot)
- \partial \phi^{(2)}(0,\cdot)\|_{L^2}
\nonumber \\
& & + C_3 T^{1/2}\bigl(\|{\tilde \phi}^{(1)}\|_{E_3(T)}
+ \|{\tilde \phi}^{(2)}\|_{E_3(T)}
+ \|\phi^{(2)}\|_{Y_3(T)}\bigr)
 \nonumber \\
& & \quad \times \bigl(
\|{\tilde \phi}^{(1)} - {\tilde \phi}^{(2)}\|_{E_1(T)\cap Y_1(T)}
+ \|\phi^{(1)} - \phi^{(2)}\|_{Y_1(T)}
\bigr).
\nonumber
\end{eqnarray}
%Here, $C_4$ is a positive constant depending on $F$, $h$, $C_{0}$, $C_{\rm S}$, and $C_2$.
\end{lemma}
\begin{prf}
If we set $\phi^* = \phi^{(1)} - \phi^{(2)}$, it satisfies
$$\partial_t^2 \phi^* - \Delta \phi^* +
 h({\tilde \phi}^{(1)})\Delta \phi^*=F(\partial {\tilde \phi}^{(1)}) - F(\partial {\tilde \phi}^{(2)})
+(h({\tilde \phi}^{(2)}) - h({\tilde \phi}^{(1)})\bigr)\Delta
\phi^{(2)}$$
By (\ref{edx12a}), we have
\begin{eqnarray*}
 \|\phi^*\|_{E_1(T)\cap Y_1(T)\cap Z_1(T)}
&\les&\|\partial \phi^*(0,\cdot)\|_{L^2}
+  \|r^{1/2}\bar\partial
h({\tilde \phi}^{(1)})\|_{L^{\infty}_{t,x}}
\|\phi^*\|_{Y_1(T)}T^{1/2}
\\
& &+ \|r^{1/4}\bigl(F(\partial {\tilde \phi}^{(1)}) - F(\partial
 {\tilde \phi}^{(2)})\bigr)\|_{L^2_{t,x}} T^{1/4}
\nonumber \\
& &+ \|r^{1/4}\bigl(h({\tilde \phi}^{(2)})
- h({\tilde \phi}^{(1)})\bigr)\Delta
\phi^{(2)}\|_{L^2_{t,x}}T^{1/4}
\end{eqnarray*}

By \eqref{ed23} and the Sobolev embedding,
we have
$$
\|r^{1/2}\bar\partial
h({\tilde \phi}^{(1)})\|_{L^{\infty}_{t,x}}\les
\|h({\tilde \phi}^{(1)})\|_{L^{\infty}_{t,x}}+ \|r^{1/2}\partial h({\tilde \phi}^{(1)})\|_{L^{\infty}_{t,x}}\les
\|{\tilde \phi}^{(1)}\|_{E_3(T)}\ ,$$
$$|\partial {\tilde \phi}^{(i)}(t,\cdot)| \les   r^{-\vep}\|\partial {\tilde \phi}^{(i)}(t,\cdot)\|_{\dot{H}^{2-\vep}}, \vep =\frac{1}{2(p-1)}\in (0, 3/2), i=1,2.$$
As $\|{\tilde \phi}^{(i)}\|_{E_3(T)} \leq C_{1}$, $p\ge 2$ and $|F(\partial {\tilde \phi}^{(1)}) - F(\partial {\tilde \phi}^{(2)})|\les |\partial {\tilde \phi}^{(1)}-\partial{\tilde \phi}^{(2)}|(|\partial {\tilde \phi}^{(1)}|^{p-1}+|\partial {\tilde \phi}^{(2)}|^{p-1})$,
we obtain
\begin{align*}
& \|r^{1/4}\bigl(F(\partial {\tilde \phi}^{(1)}) - F(\partial
 {\tilde \phi}^{(2)})\bigr)\|_{L^2_{t,x}} \\
\les& \|r^{1/4-\vep(p-1)}(\partial {\tilde \phi}^{(1)}-\partial{\tilde \phi}^{(2)})\|_{L^{2}(S_{T})}
(\|{\tilde \phi}^{(1)}\|_{E_3(T)}^{p-1}+\|{\tilde \phi}^{(2)}\|_{E_3(T)}^{p-1})\\
\les & T^{1/4}\|{\tilde \phi}^{(1)}-{\tilde \phi}^{(2)}\|_{Y_1(T)}(\|{\tilde \phi}^{(1)}\|_{E_3(T)}+\|{\tilde \phi}^{(2)}\|_{E_3(T)}).
\end{align*}

Since
\begin{equation*}
|h({\tilde \phi}^{(2)})-h({\tilde \phi}^{(1)})|\les
|{\tilde \phi}^{(2)}-{\tilde \phi}^{(1)}|
\les r^{-1}\|{\tilde \phi}^{(2)}-{\tilde \phi}^{(1)}\|_{{\dot{H}}^{1}_{{\rm rad}}({\mathbb{R}}^{4})}\ ,
\end{equation*}
we have
\begin{eqnarray*}
\|r^{\frac{1}{4}}(h({\tilde \phi}^{(2)})-h({\tilde \phi}^{(1)}))\Delta\phi^{(2)}\|_{L^{2}(S_{T})}
&\les & \|r^{-\frac{3}{4}}\Delta\phi^{(2)}\|_{L^{2}(S_{T})}\|{\tilde \phi}^{(2)}-{\tilde \phi}^{(1)}\|_{E_1(T)}\\
&\les&T^{1/4}\|\phi^{(2)}\|_{Y_3(T)}\|{\tilde \phi}^{(2)}-{\tilde \phi}^{(1)}\|_{E_1(T)}\ ,
\end{eqnarray*} where, in the second inequality we used
$$
\|r^{-\frac{3}{4}}\Delta\phi^{(2)}\|_{L^{2}(S_{T})}
\leq  \|r^{-\frac{1}{4}}\Delta\phi^{(2)}\|_{L^{2}(S_{T})}^{\frac{1}{2}}
\|r^{-\frac{5}{4}}\Delta\phi^{(2)}\|_{L^{2}(S_{T})}^{\frac{1}{2}}
\leq  T^{1/4}\|\phi^{(2)}\|_{Y_3(T)}.
$$
This completes the proof of Lemma \ref{ld9}.
\end{prf}
\begin{lemma}
\label{edx123}
Let $(f,g) \in H^3_{\rm rad}({\mathbb R}^4)\times H^2_{\rm rad}({\mathbb R}^4)$.
There exist positive constants $\varepsilon_2, A_2$, so that
the sequence $\{\phi_k\}$ defined in Lemma {\rm \ref{exd11111}}
converges to a function $\phi \in X_1(T) \cap Y_1(T)\cap Z_1(T)$
if $\|\nabla f\|_{H^2} + \|g\|_{H^2} = \varepsilon \leq \varepsilon_2$
and $T \leq A_2 \varepsilon^{-2}$.
In addition, $\phi \in L^{\infty}([0,T]; H^3({\mathbb
 R}^4)) \cap C^{0,1}([0,T]; H^2({\mathbb
 R}^4)) \cap Y_3(T)\cap Z_3(T)$
and it is a unique weak solution for
the initial value problem \eqref{ea1}-\eqref{ea2}.
\end{lemma}
\begin{prf}
Let $\varepsilon \leq \varepsilon_{1}$ and $T < \infty$,
by Lemma {\rm \ref{exd11111}},
${\phi_{k}}$ is well defined and
\begin{equation*}
\|\phi_k\|_{E_3(T)}
 + \|\phi_k\|_{Y_3(T)}
 + \|\phi_k\|_{Z_3(T)}
\leq C_2 \varepsilon.
\end{equation*}
If
\begin{equation*}
2C_{2} C_{3} \varepsilon T^{\frac{1}{2}}\leq \frac{1}{4},
\end{equation*}
then by Lemma \ref{ld9}, we have for $k\ge 2$,
\begin{eqnarray}
\lefteqn{\|\phi_k - \phi_{k-1}\|_{E_1(T)}
+ \|\phi_k - \phi_{k-1}\|_{Y_1(T)\cap Z_1(T)}
} \nonumber \\
&\leq& \frac{4}{3}C_3\|\partial \phi_k(0,\cdot)
- \partial \phi_{k-1}(0,\cdot)\|_{L^2}
 \nonumber \\
& & + \frac{1}{3}\bigl(\|\phi_{k-1} - \phi_{k-2}\|_{E_1(T)}
 + \|\phi_{k-1} - \phi_{k-2}\|_{Y_1(T)}
\bigr),
\nonumber
\end{eqnarray}
which, together with (\ref{ed12}), implies
\begin{equation}
 \sum_{k=1}^{\infty}\bigl(\|\phi_k - \phi_{k-1}\|_{E_1(T)}
+ \|\phi_k - \phi_{k-1}\|_{Y_1(T)\cap Z_1(T)}\bigr) < \infty.
\nonumber
\end{equation}
Hence, $\{\phi_k\}$ is a Cauchy sequence in $X_1(T)\cap Y_1(T)$
as long as $T \leq A_2 \varepsilon^{-2}$ and $\varepsilon \leq \varepsilon_2$,
where $A_2$ and $\varepsilon_2$ are chosen so that
\begin{equation}
 \label{ed54}
 \varepsilon_2 \leq \varepsilon_1,\quad
\varepsilon_2 \leq A_2^{1/2},
\quad A_2 = \frac{1}{64 C_2^2 C_3^2}.
\end{equation}
Now it remains to prove that its limit $\phi$ belongs to
$L^{\infty}([0,T]; H^3) \cap C^{0,1}([0,T]; H^2) \cap Y_3(T)\cap Z_3(T)$. Since $\{\partial_{\alpha} \phi_k(t,\cdot)\}$ is bounded in
$H^2({\mathbb R}^4)$ and $\{\partial_{\alpha} \phi_{k}(t,\cdot)\}$ converges to
$\partial_{\alpha} \phi(t,\cdot)$ in $L^2({\mathbb R}^4)$, we see that
$\{\partial_{\alpha} \phi_{k}(t,\cdot)\}$ has a unique limit point
$\partial_{\alpha} \phi(t,\cdot)$ with respect to the weak topology of $H^2({\mathbb R}^4)$. Therefore,
\begin{equation}
 \label{ed71}
 \|\partial \phi(t,\cdot)\|_{H^2({\mathbb R}^4)} \leq \liminf_{k \to \infty}
 \|\partial \phi_{k}(t,\cdot)\|_{H^2({\mathbb R}^4)} \leq C_2 \varepsilon
\end{equation}
for $0 \leq t \leq T \leq A_2\varepsilon^{-2}$,
which shows that $\phi \in L^{\infty}([0,T]; H^3)
\cap C^{0,1}([0,T]; H^2)$.
By similar arguments, we conclude that
$\phi \in Y_3(T)\cap Z_3(T)$ and
\begin{equation}
 \|\phi\|_{Y_3(T)} + \|\phi\|_{Z_3(T)} \leq \liminf_{k \to \infty}
\bigl(\|\phi_k\|_{Y_3(T)} +  \|\phi_k\|_{Z_3(T)}\bigr) \leq M_1 \varepsilon.
\nonumber
\end{equation}
By definition of $\phi_k$, it is clear that $\phi$ solves
\eqref{ea1}-\eqref{ea2}.
For uniqueness, we need only to apply
Lemma \ref{ld9} again.
This completes the proof of Lemma \ref{edx123}.
\end{prf}

\subsection{Global existence}
Now we are ready to complete the proof of Theorem \ref{ta1}.
So far, we have proven that there exists a unique weak solution for
the initial value problem (\ref{ea1}), (\ref{ea2}) for $T \leq
A_2\varepsilon^{-2}$,
where $\varepsilon = \|\nabla f\|_{H^2} + \|g\|_{H^2}$ and $\varepsilon
\leq \varepsilon_2$.
It remains to show that the solution actually extends
to any finite time $T < \infty$.

Let $T_1 = A_2\varepsilon_2^{-2}$,
we define a sequence $\{\psi^{T_1}_k\}$ by \eqref{ed14},
with the initial data $(f_k, g_k)$ replaced by
$$
 \tilde f_k = \rho_{k}*\phi(T_1,\cdot),\quad
 \tilde g_k = \rho_{k}*\partial_t\phi(T_1,\cdot)\ .
$$
Note that $\|\nabla \phi(T_1,\cdot)\|_{H^2}
+ \|\partial_t \phi(T_1,\cdot)\|_{H^2} \leq C_2\varepsilon$.
Let $\varepsilon$ satisfies $C_2\varepsilon \leq \varepsilon_2$ and $A_2(C_2\varepsilon)^{-2}\ge 1$,
then by Lemma \ref{exd11111} and Lemma \ref{edx123}, we see that,
with $T =A_2(C_2\varepsilon)^{-2}$, $$\|\psi_k^{T_1}\|_{E_3(T)} + \|\psi_k^{T_1}\|_{Y_3(T)}
+ \|\psi_k^{T_1}\|_{Z_3(T)} \leq C_2^2\varepsilon\ ,$$
and $\{\psi^{T_1}_k\}$ converges to a function $\psi^{T_1}$ in
$X_1(T) \cap Y_1(T) \cap Z_1(T)$.
In addition, $\psi^{T_1} \in X_3(T) \cap Y_3(T) \cap Z_3(T)$ and
$\|\psi^{T_1}\|_{E_3(T)\cap Y_3(T) \cap Z_3(T)}
 \leq C_2^2\varepsilon$.

 On the other hand, if
we define another sequence $\{\phi_k^{T_1}\}$ on $[0,T]\times {\mathbb R}^4$,
defined by $\phi^{T_1}_k(t,x) = \phi_k(t+T_1,x)$.
Then it follows from Lemma \ref{ld9} that for $T\ge 1$,
\begin{eqnarray*}
\|\psi_k^{T_1} - \phi_k^{T_1}\|_{E_1(T)\cap Y_1(T)}
&\leq& C_3\|\partial \psi^{T_1}_k(0,\cdot)
- \partial \phi^{T_1}_k(0,\cdot)\|_{L^2} \nonumber \\
& & + C_3T^{1/2} (\|\psi^{T_1}_{k-1}\|_{E_3(T)}
+ \|\phi^{T_1}_{k-1}\|_{E_3(T)}
+ \|\phi^{T_1}_{k}\|_{Y_3(T)}) \\
& & \times \bigl(
\|\psi_{k-1}^{T_1} - \phi_{k-1}^{T_1}\|_{E_1(T)\cap Y_1(T)}
+ \|\psi_k^{T_1} - \phi_k^{T_1}\|_{Y_1(T)}
\bigr).
\end{eqnarray*}
Taking $A_3>0$ and $\varepsilon_0>0$ so that
$$
C_2 \varepsilon_0 \leq \varepsilon_2,
 \varepsilon_0\leq A_3^{1/2},
 A_3 \leq A_2C_2^{-2},
C_3(2C_2 + C_2^2)\sqrt{A_3}\leq \frac{1}{4}. $$
Then for $T_2 = A_3\varepsilon^{-2}_0$($\in [1,T]$) and $\varepsilon \leq \varepsilon_0$, we have
\begin{equation}
 C_3 (2C_2 \varepsilon + C_2^2 \varepsilon) T_2^{1/2} \leq \frac{1}{4},
\nonumber
\end{equation}
and thus
$$
\|\psi^{T_1}_k - \phi^{T_1}_k\|_{E_1(T_2)\cap Y_1(T_2)}
\le \frac{4}{3}C_3\|\partial \psi^{T_1}_k(0,\cdot)
- \partial \phi^{T_1}_k(0,\cdot)\|_{L^2}
+ \frac{1}{3} \|\psi^{T_1}_{k-1} - \phi^{T_1}_{k-1}\|_{E_1(T_2)\cap Y_1(T_2)}$$

Hence we have $\|\psi^{T_1}_{k} - \phi^{T_1}_{k}\|_{E_1(T_2)\cap Y_1(T_2)} \to 0$ as $k \to
\infty$, which implies $\phi^{T_1}_k \to \psi^{T_1}$
in $X_1(T_2) \cap Y_1(T_2)$.
Therefore, we have shown that
there exists an extension of $\phi \in X_1(T_1) \cap
Y_1(T_1)$ such that
$\phi_k \to \phi$ in $X_1(T_1 + T_2)
\cap Y_1(T_1 + T_2)$, while we still have the uniform bound
$\|\phi\|_{E_3(T_1 + T_2)\cap Y_3(T_1 + T_2)\cap Z_3(T_1 + T_2)} \leq C_2\varepsilon$.
Repeating this argument, we can extend the solution to any $T < \infty$.
This completes the proof of Theorem \ref{ta1}.

\bibliographystyle{plain1}

\end{document}